\documentclass{article}
\usepackage{amsmath,amssymb,amsthm}
\begin{document}
\title{\bf Interval valued intuitionistic $(S,T)$-fuzzy  $H_v$-submodules}
\normalsize \author{Jianming Zhan \\ {\small Department of
Mathematics,
  Hubei Institute for Nationalities,}\\ {\small Enshi, Hubei Province,
   445000, P. R. China}\\ {\small Email:
zhanjianming@hotmail.com}
\\[10pt]
Wieslaw A. Dudek \\ {\small Institute of Mathematics and Computer
Science,}\\ {\small Wroclaw University of Technology,}\\ {\small
Wyb. Wyspianskiego 27, 50-370 Wroclaw, Poland}\\ {\small Email:
dudek@im.pwr.wroc.pl}}
 \date{}\maketitle

{\bf Abstract.} On the basis of the concept of the interval valued
intuitionistic fuzzy sets introduced by K.Atanassov, the notion of
interval valued intuitionistic fuzzy $H_v$-submodules of an
$H_v$-module with respect to $t$-norm $T$ and $s$-norm $S$ is
given and the characteristic properties are described. The
homomorphic image and the inverse image are investigated.In
particular, the connections between interval valued intuitionistic
$(S,T)$-fuzzy $H_v$-submodules and interval valued intuitionistic
$(S,T)$-fuzzy submodules are discussed.



\maketitle

\section{Introduction }

The concept of hyperstructure was introduced in 1934 by Marty
\cite{5} at the 8th Congress of Scandinavian Mathematicians.
Hypersructures have many applications to several branches of both
pure and applied sciences (see for example \cite{21} and
\cite{22}). Vougiouklis \cite{6, 7} introduced a new class of
hyperstructures called now $H_v$-structures, and Davvaz \cite{4}
surveyed the theory of $H_v$-structures. After the introduction of
fuzzy sets by Zadeh \cite{9}, there have been a number of
generalizations of this fundamental concept. The notion of
intuitionistic fuzzy sets introduced by Atanassov \cite{1} is one
among them. For more details on intuitionistic fuzzy sets, we
refer the reader to \cite{1, 2, 20}. In 1975, Zadeh \cite{10}
introduced the concept of interval valued fuzzy subsets, where the
values of the membership functions are intervals of numbers
instead of the numbers.

Such fuzzy sets have some applications in the technological scheme
of the functioning of a silo-farm with pneumatic transportation in
a plastic products company and in medicine (see the book
\cite{20}).

In this paper, we introduce the notion of interval valued
intuitionistic $(S,T)$-fuzzy $H_v$-submodules of an $H_v$-module
and describe the characteristic properties. We give the
homomorphic image and the inverse image. In particular, we discuss
the connections between interval valued intuitionistic
$(S,T)$-fuzzy $H_v$-submodules and interval valued intuitionistic
$(S,T)$-fuzzy submodules.

\section{Preliminaries}

In this section, we recall some basic definitions for the sake of
completeness.

As it is well known \cite{6}, a {\it hyperstructure} is a
non-empty set $H$ together with a map $\cdot: H\times H\rightarrow
P^*(H)$, called a {\it hyperoperation}, where $P^*(H)$ is the
family of all non-empty subsets of $H$. The image of pair $(x,y)$
is denoted by $x\cdot y$. If $x\in H$ and $A,B\subseteq H$, then
by $A\cdot B$, $A\cdot x$ and $x\cdot B$ we mean
\[
A\cdot B=\displaystyle\bigcup_{a\in A, b\in B}\!\!\! a\cdot b, \ \
A\cdot x=A\cdot \{ x\} \ \ {\rm and } \ \ x\cdot B=\{x\}\cdot B,
\]
respectively.

\paragraph{Definition 2.1.} A hyperstructure $(H,\cdot )$ is
called an {\it $H_v$-semigroup } if
\[
(x\cdot (y\cdot z))\cap ((x\cdot y)\cdot z)\neq\emptyset \ \ {\rm
for \ all } \ x,y,z\in H.
\]
An $H_v$-semigroup in which \ $a\cdot H=H\cdot a=H$ is valid for
every $a\in H$ is called an {\it $H_v$-group}.

The last condition means that for any $a,h\in H$ there exist
$u,v\in H$ such that $h\in a\cdot u$ and $h\in v\cdot a$. An
$H_v$-group $(H,\cdot )$ satisfying for all $x,y\in H$ the
condition $x\cdot y\cap y\cdot x\neq\emptyset$ is called {\it weak
commutative}.

\paragraph{Definition 2.2.} An {\it $H_v$-ring} is a system $(R,+,\cdot)$ with
two hyperoperations satisfying the following axioms:

(i) \ $(R,+)$ is an $H_v$-group;,

(ii) \ $(R,\cdot)$ is an $H_v$-semigroup,

(iii) \ the multiplication $\cdot$ is weak distributive with
respect to the addition $+$, i.e.,
\[
\begin{array}{l}
( x\cdot (y+z))\cap (x\cdot y+x\cdot z)\neq\emptyset ,\\[2pt]
( (x+y)\cdot z)\cap (x\cdot z+y\cdot z)\neq\emptyset
\end{array}
\]
for all $x,y,z\in R$.

\paragraph{Definition 2.3 (\cite{11}).}
A non-empty set $M$ is an {\it $H_v$-module over an $H_v$-ring}
$R$ if $(M,+)$ is a weak commutative $H_v$-group and there exists
the map \ $\cdot : R\times M\rightarrow P^*(M)$, \ $(r,x)\mapsto
r\cdot x$, such that for all $a\in R$ and $x,y\in M$, we have

(i) \ $(a\cdot (x+y))\cap (a\cdot x+a\cdot y)\ne\emptyset$,

(ii) \ $((a+b)\cdot x)\cap (a\cdot x+b\cdot x)\ne\emptyset$,

(iii) \ $((a\cdot b)\cdot x)\cap (a\cdot(b\cdot x))\ne\emptyset$.

A non-empty subset $S$ of $M$ is an {\it $H_v$-submodule } of $M$
if $(S,+)$ is an $H_v$-subgroup of $(M,+)$ and $R\cdot S\subseteq
S$.

It is clear that an arbitrary ring (module) will be an $H_v$-ring
($H_v$-module) if we identify $x$ with $\{x\}$. Others interesting
examples are given below.

\paragraph{Example 2.4.} Let $(M,+,\cdot )$ be an ordinary module over a ring
$R$ with a center $Z(R)$. On $R\times M$ we can define three
hyperoperations $P^*$, $P_+$ and $P_+^*$ putting for all $(r,x)\in
R\times M$:
\begin{itemize}
\item[(1)] \ $rP^*x=(rP)x$ \ if \ $P\subseteq R$,
\item[(2)] \ $rP_+x=r(P+x)$ \ if \ $P\subseteq M$,
\item[(3)] \ $rP^*_+x=(rP_1)(P_2+x)$ \ if \ $P_1\subseteq R$\ and\ $P_2 \subseteq M$.
\end{itemize}
Then, as it is not difficult to verify,
\begin{itemize}
\item[(a)] \ $(M,+,P^*)$ is an $H_v$-module over $R$, if there exists $p\in P\cap Z(R)$ such that $p^2\in P$,
\item[(b)] \ $(M,+,P_+)$ is an $H_v$-module over $R$, if the zero $0$ of $(M,+)$ belongs to
$P\subseteq M$,
\item[(c)] \ $(M,+,P^*_+)$ is an $H_v$-module over $R$, if there exist
$p_1\in P_1\cap Z(R)$ such that $p^2_1=p_1$ and $p_2\in
P_2\subseteq M$ such that $p_1\cdot p_2=0$.
\end{itemize}

\medskip

According to Zadeh \cite{9}, a {\it fuzzy set} $\mu_F$ defined on
a non-empty set $X$, i.e. a map $\mu_F : X\rightarrow [0,1]$, can
be identified with the set $F=\{(x,\mu_F(x))\ |\ x\in X\}$.

\paragraph{Definition 2.5 (\cite{4}).} A fuzzy set $F$ of an $H_v$-module $M$
over an $H_v$-ring $R$ is said to be a {\it fuzzy $H_v$-submodule}
of $M$ if:
\begin{enumerate}
\item[(i)] \ $\min\{\mu_F(x),\mu_F(y)\}\le \inf\limits_{\alpha\in
x+y}\mu_F(\alpha)$ \ \ for all $x,y\in M$,

\item[(ii)] \ for all $x,a\in M,$ there exists $y\in M$ such that
$x\in a+y$ and

\ $\min\{\mu_F(a),\mu_F(x)\}\le\mu_F(y)$,

\item[(iii)] \ for all $x,a\in M,$ there exists $z\in M$ such that
$x\in z+a$ and

\ $\min\{\mu_F(a),\mu_F(x)\}\le\mu_F(z)$,

\item[(iv)] \ $\mu_F(x)\le\inf\limits_{\alpha\in r\cdot
x}\mu_F(\alpha)$ \ \ for all $r\in R$ and $x\in M$.
\end{enumerate}

\medskip

By an {\it interval number } $\tilde{a}$ we mean (cf. \cite{2}) an
interval $[a^-,a^+]$, where $0\le a^-\le a^+ \le 1$. The set of
all interval numbers is denoted by $D[0,1]$. The interval $[a,a]$
is identified with the number $a\in [0,1]$.

For interval numbers \ $\widetilde{a}_i=[a_i^-,a_i^+]\in D[0,1]$,
$i\in I,$ we define
 \[
\inf\widetilde{a}_i=[\bigwedge\limits_{i\in I} a_i^-,
\bigwedge\limits_{i\in I} a_i^+], \ \ \ \ \ \
\sup\widetilde{a}_i=[\bigvee\limits_{i\in I} a_i^-,
\bigvee\limits_{i\in I} a_i^+]
 \]
and put
\begin{enumerate}
\item[(1)] \ $\widetilde{a}_1\le \widetilde{a}_2\Longleftrightarrow
a_1^-\le a_2^-$ and $a_1^+\le a_2^+$,

\item[(2)] \ $\widetilde{a}_1=\widetilde{a}_2\Longleftrightarrow  a_1^-=a_2^-$ and
$a_1^+=a_2^+$,

\item[(3)] \ $\widetilde{a}_1<\widetilde{a}_2\Longleftrightarrow
\widetilde{a}_1\le \widetilde{a}_2$ and $\widetilde{a}_1\ne
\widetilde{a}_2$,

\item[(4)] \ $k\widetilde{a}=[ka^-,ka^+]$, whenever $ 0\le k\le 1$.
\end{enumerate}
In is clear that $(D[0,1],\le,\vee,\wedge)$ is a complete lattice
with $0=[0,0]$ as the least element and $1=[1,1]$ as the greatest
element.

\medskip

By an {\it interval valued fuzzy set} $F$ on $X$ we mean (sf.
\cite{10}) the set
 \[
F=\{(x, [\mu^-_F(x),\mu^+_F(x)])\,|\,x\in X\},
\]
where $\mu^-_F $ and $\mu^+_F$ are two fuzzy subsets of $X$ such
that $\mu^-_F(x)\le \mu^+_F(x)$ for all $x\in X.$ Putting
$\mu_F(x)=[\mu^-_F(x),\mu^+_F(x)]$, we see that
$F=\{(x,\mu_F(x))\,|\,x\in X\}$, where $\mu_F: X\rightarrow
D[0,1].$

\medskip

As it is well-known, any function $\delta: [0,1]\times
[0,1]\rightarrow [0,1]$ such that $\delta(x,y)=\delta(y,x)$,
$\delta(x,x)=x$, $\delta(\delta(x,y),z)=\delta(x,\delta(y,z))$ and
$\delta(x,u)\le \delta(x,w)$ for all $x,y,z,u,w\in [0,1]$, where
$u\le w$ is called an {\it idempotent $t$-norm} if
$\delta(x,1)=x$, and an {\it idempotent $s$-norm} if
$\delta(1,1)=1$ and $\delta(x,0)=x$ for all $x\in [0,1]$.

If $\delta$ is an idempotent $t$-norm ($s$-norm), then the mapping
$\Delta: D[0,1]\times D[0,1]\rightarrow D[0,1]$ defined by
$\Delta(\widetilde{a}_1,\widetilde{a}_2)=[\delta(a_1^-,a_2^-),
\delta(a_1^+,a_2^+)]$ is, as it is not difficult to verify, an
idempo\-tent $t$-norm ($s$-norm, respectively) and is called an
{\it idempotent interval $t$-norm} ({\it $s$-norm}, respectively).

\medskip

According to Atanassov (cf. \cite{1, 20}) an {\it interval valued
intuitionistic fuzzy set} on $X$ is defined as the object of the
form
$$
A=\{(x, \widetilde{M}_A(x), \widetilde{N}_A (x))\mid x\in X\},
$$
where $\widetilde{M}_A(x)$ and $\widetilde{N}_A(x)$ are interval
valued fuzzy sets on $X$ such that
$$
0\le\sup\widetilde{M}_A(x)+\sup\widetilde{N}_A(x)\le 1 \ \ \ \
{\rm for\; all} \ \ x\in X .
$$
For the sake of simplicity, in the sequel such interval valued
intuitionistic fuzzy set will be denoted by
$A=(\widetilde{M}_A,\widetilde{N}_A)$.

\section{Interval valued intuitionistic $(S,T)$-fuzzy $H_v$-submodules}

In what follow, let $M$ denote an $H_v$-module over an $H_v$-ring
$R$ unless otherwise specified.

\paragraph{Definition 3.1.}
Let $T$ (resp. $S$) be an idempotent interval $t$-norm (resp.
$s$-norm). An interval valued intuitionistic fuzzy set
$A=(\widetilde{M}_A,\widetilde{N}_A)$ of $M$ is called an {\it
interval valued intuitionistic $(S,T)$-fuzzy $H_v$-submodule} of
$M$ if the following condition hold:
\begin{enumerate}
\item[(1)] \ $T(\widetilde{M}_A(x),\widetilde{M}_A(y))\le
\inf\limits_{\alpha\in x+y} \widetilde{M}_A(\alpha)$  and
$S(\widetilde{N}_A(x),\widetilde{N}_A(y))\ge
\sup\limits_{\alpha\in x+y} \widetilde{N}_A(\alpha), \ \forall
x,y\in M$,

\item[(2)] \ for all $x,a\in M$ there exists $y\in M$ such that
$x\in a+y$, $T(\widetilde{M}_A(x),\widetilde{M}_A(a))\le
\widetilde{M}_A(y)$ and
$S(\widetilde{N}_A(x),\widetilde{N}_A(a))\ge\widetilde{N}_A(y)$,

\item[(3)] \ for all $x,a\in M$ there exists $z\in M$ such that
$x\in z+a$, $T(\widetilde{M}_A(x),\widetilde{M}_A(a))\le
\widetilde{M}_A(z)$ and
$S(\widetilde{N}_A(x),\widetilde{N}_A(a))\ge\widetilde{N}_A(z)$,

\item[(4)] \ $\widetilde{M}_A(x)\le\inf\limits_{\alpha\in r\cdot x}
\widetilde{M}_A(\alpha)$ and $\widetilde{N}_A(x)\ge
\sup\limits_{\alpha\in r\cdot x}\hat {N}_A(\alpha)$ \ for all
$x\in M$ and $r\in R$.
\end{enumerate}
With any interval valued intuitionistic fuzzy set
$A=(\widetilde{M}_A,\widetilde{N}_A)$ of $M$ are connected two
levels:
\[
U(\widetilde{M}_A;[t,s])=\{x\in X\,|\,\widetilde{M}_A(x)\ge
[t,s]\},
\]
and
\[
L(\widetilde{N}_A;[t,s])=\{x\in X\,|\,\widetilde{N}_A(x)\le
[t,s]\}.
\]
\paragraph{Theorem 3.2.}
{\em Let $T$ $($resp. $S\,)$ be an idempotent interval $t$-norm
$($resp. $s$-norm$)$, then $A=(\widetilde{M}_A,\widetilde{N}_A)$
is an interval valued intuitionistic $(S,T)$-fuzzy $H_v$-submodule
of $M$ if and only if every for all $t,s\in [0,1]$, $t\le s$,
$U(\widetilde{M}_A;[t,s])$ and $L(\widetilde{N}_A;[t,s])$ are
$H_v$-submodules of $M$.}

 \medskip
\noindent{\it Proof}. Let $A=(\widetilde{M}_A,\widetilde{N}_A)$ be
an interval valued intuitionistic $(S,T)$-fuzzy $H_v$-submodule of
$M.$ Then for every $x,y\in U(\widetilde{M}_A;[t,s])$ we have
$\widetilde{M}_A(x)\ge [t,s]$  and $\widetilde{M}_A(y)\ge [t,s]$.
Hence $T(\widetilde{M}_A(x),\widetilde{M}_A(y))\ge
T([t,s],[t,s])=[t,s]$, and so $\inf\limits_{\alpha\in x+y}
\widetilde{M}_A(\alpha)\ge [t,s]$. Therefore $\alpha\in
U(\widetilde{M}_A;[t,s])$ for every $\alpha\in x+y$, so
$x+y\subseteq U(\widetilde{M}_A;[t,s])$. Thus, for every $a\in
U(\widetilde{M}_A;[t,s])$, we have $a+
U(\widetilde{M}_A;[t,s])\subseteq  U(\widetilde{M}_A;[t,s])$. On
the other hand, for $x,a\in U(\widetilde{M}_A;[t,s])$ there exists
$y\in H$ such that $x\in a+ y$ and
$T(\widetilde{M}_A(x),\widetilde{M}_A(a))\le \widetilde{M}_A(y)$.
But $T(\widetilde{M}_A(x),\widetilde{M}_A(a))\ge [t,s]$ for all
$x,a\in U(\widetilde{M}_A;[t,s])$, so
$\widetilde{M}_A(y)\ge[t,s]$, that is, $y\in
U(\widetilde{M}_A;[t,s])$. Whence
$U(\widetilde{M}_A;[t,s])\subseteq a+U(\widetilde{M}_A;[t,s])$,
and, in the consequence $U(\widetilde{M}_A;[t,s])=
a+U(\widetilde{M}_A;[t,s])$. Similarly, we can prove that
$U(\widetilde{M}_A;[t,s])= U(\widetilde{M}_A;[t,s])+a$. That is,
$a+U(\widetilde{M}_A;[t,s])= U(\widetilde{M}_A;[t,s])=
U(\widetilde{M}_A;[t,s])+a$. This proves that
$(U(\widetilde{M}_A;[t,s]),+)$ is an $H_v$-subgroup of $(M,+)$.

If $r\in R$ and $x\in U(\widetilde{M}_A;[t,s])$, then
$\widetilde{M}_A(x)\ge [t,s]$, which means that
$\inf\limits_{\alpha\in r\cdot x}\widetilde{M}_A(\alpha)\ge
[t,s]$. So, $\alpha\in U(\widetilde{M}_A;[t,s])$ for every
$\alpha\in r\cdot x$. Therefore, $r\cdot x\subseteq
U(\widetilde{M}_A;[t,s])$, i.e. $r\cdot
U(\widetilde{M}_A;[t,s])\subseteq U(\widetilde{M}_A;[t,s])$. This
proves that $U(\widetilde{M}_A;[t,s])$ is an $H_v$-submodule of
$M$. Similarly, we can show that $L (\widetilde{N}_A;[t,s])$ is an
$H_v$-submodule of $M$.

Conversely, assume that for every $[t,s]\in D[0,1]$ any non-empty
$U(\widetilde{M}_A;[t,s])$ is an $H_v$-submodule of $M.$ If
$[t_0,s_0]=T(\widetilde{M}_A(x),\widetilde{M}_A(y))$ for some
$x,y\in H$, then $x,y\in U(\widetilde{M}_A;[t_0,s_0])$, and so
$x+y\subseteq U(\widetilde{M}_A;[t_0,s_0])$. Therefore $\alpha\in
U(\widetilde{M}_A;[t_0,s_0])$ for every $\alpha\in x+y$, and so
$\inf\limits_{\alpha\in x+y} \widetilde{M}_A(\alpha)\ge
T(\widetilde{M}_A(x),\widetilde{M}_A(y))$. Now, if
$[t_1,s_1]=T(\widetilde{M}_A(a),\widetilde{M}_A(x))$ for some
$a,x\in H$, then $a+x\in U(\widetilde{M}_A;[t_1,s_1]) $, so there
exists $y\in U(\widetilde{M}_A;[t_1,s_1])$ such that $x\in a+y$.
But for $y\in U(\widetilde{M}_A;[t_1,s_1])$ we have
$\widetilde{M}_A(y)\ge [t_1,s_1]$, whence $\widetilde{M}_A(y)\ge
T(\widetilde{M}_A(a),\widetilde{M}_A(x))$. Similarly, we can show
that for $a,x\in H$ there exists $z\in H$ such that $x\in z+a$ and
$\widetilde{M}_A(z)\ge T(\widetilde{M}_A(a),\widetilde{M}_A(x))$.
If $[t_2,s_2]=\widetilde{M}_A(x)$ for some $x\in M$, then $x\in
U(\widetilde{M}_A;[t_2,s_2]),$ and so $r\cdot x\subseteq
U(\widetilde{M}_A;[t_2,s_2])$ for every $r\in R$. Therefore for
every $\alpha\in r\cdot x$, we have $\alpha\in
U(\widetilde{M}_A;[t_2,s_2]),$ consequently
$\inf\limits_{\alpha\in r\cdot x}\widetilde{M}_A(\alpha)\ge
[t_2,s_2]=\widetilde{M}_A(x)$. This proves that $\widetilde{M}_A$
is an interval valued $T$-fuzzy $H_v$-submodule of $M.$

Similarly, we can show that $\widetilde{N}_A$ is an interval
valued $S$-fuzzy $H_v$-submodule of $M.$ Therefore,
$A=(\widetilde{M}_A,\widetilde{N}_A)$ is an interval valued
intuitionistic $(S,T)$-fuzzy $H_v$-submodule of $M$.
\hfill$\Box{}$

\paragraph{Definition 3.3.} Let $f: X\rightarrow Y$ be a mapping and
$A=(\widetilde{M}_A,\widetilde{N}_A)$,
$B=(\widetilde{M}_B,\widetilde{N}_B)$ an interval valued
intuitionistic sets $X$ and $Y,$ respectively. Then  the {\it
image} $f[A]=(f(\widetilde{M}_A),f(\widetilde{N}_A))$ of $A$ is
the interval valued intuitionistic fuzzy set of $Y$ defined by
\[
f(\widetilde{M}_A)(y)= \left\{\begin{array}{clll}\sup\limits_{z\in
f^{-1}(y)}\widetilde{M}_A(z) &&
\mbox{ if\ \ } f^{-1}(y)\ne\emptyset\\[12pt] [0,0] &&
 \mbox{ otherwise}
 \end{array}\right.
\]
\[
f(\widetilde{N}_A)(y)=\left\{\begin{array}{clll}\inf\limits_{z\in
f^{-1}(y)} \widetilde{N}_A(z) && \mbox{ if\ \ }
f^{-1}(y)\ne\emptyset\\[12pt] [1,1] &&
\mbox{ otherwise}\end{array}\right.
\]
for all $y\in Y.$

The inverse image $f^{-1}(B) $ of $B$ is an interval valued
intuitionistic fuzzy set defined by
$f^{-1}(\widetilde{M}_B(x)=\widetilde{M}_B(f(x)),$ \
$f^{-1}(\widetilde{N}_B)(x)=\widetilde{N}_B(f(x))$ for all $x\in
X$.

\paragraph{Definition 3.4 ([3]).}\rm Let $M$ and $N$ be two $H_v$-modules
over an $H_v$-ring $R$. A mapping $f: M\rightarrow N$ is called an
{\it $H_v$-homomorphism} or {\it weak homomorphism} if for all
$x,y\in M$ and $\ r\in R$, the following relations hold: $\
f(x+y)\bigcap(f(x)+f(y))\ne\emptyset$ \ and \ $f(r\cdot x)\bigcap
r\cdot f(x)\ne\emptyset$. $f$ is called an {\it inclusion
homomorphism} if $f(x+y)\subseteq f(x)+f(y)$ and $f(r\cdot
x)\subseteq r\cdot f(x)$ for all $x,y\in M $ and $r\in R$.
Finally, $f$ is called a {\it strong homomorphism} if for all
$x,y\in M$ and $x\in R$, we have $f(x+y)=f(x)+f(y)$ and $f(r\cdot
x)=r\cdot f(x)$.

\paragraph{Lemma 3.5 ([3]).} {\em Let $M_1$ and $M_2$ be two
$H_v$-modules over an $H_v$-ring $R$ and $f: M_1\rightarrow M_2$ a
strong epimorphism. If $\,N$ is an $H_v$-submodule of $M_2$, then
$f^{-1}(N)$ is an $H_v$-submodule of $\,M_1$.}

\paragraph{Theorem 3.6.} {\em
Let $M_1$ and $M_2$ be two $H_v$-modules, $f$ a strong
homomorphism from $H_1$ into $H_2$ and $T$ $($resp. $S\,)$ an
idempotent interval $t$-norm $($resp. $s$-norm$)$.
\begin{enumerate}
\item[$(i)$] If $A=(\widetilde{M}_A,\widetilde{N}_A)$ is an interval
valued intuitionistic $(S,T)$-fuzzy $H_v$-submodule of $M_1$, then
the image $f[A]$ of $A$ is an interval intuitionistic
$(S,T)$-fuzzy $H_v$-submodule of $M_2$.

\item[$(ii)$] If $B=(\widetilde{M}_B,\widetilde{N}_B)$ be an interval
valued intuitionistic $(S,T)$-fuzzy $H_v$-submodule of $M_2$, then
the inverse image $f^{-1}(B)$ of $B$ is an interval valued
intuitionistic $(S,T)$-fuzzy $H_v$-submodule of $M_1$.
\end{enumerate} }

\noindent
 {\it Proof.}\
$(i)$ Let $A=(\widetilde{M}_A,\widetilde{N}_A)$ be an interval
valued intuitionistic $(S,T)$-fuzzy $H_v$-submodule of $M_1$. By
Theorem 3.2, $U(\widetilde{M}_A;[t,s])$ and
$L(\widetilde{N}_A;[t,s])$ are $H_v$-submodules of $M_1$ for every
$[t,s]\in D[0,1]$. Therefore, by Lemma 3.5,
$f(U(\widetilde{M}_A;[t,s])$ and $f(L(\widetilde{N}_A;[t,s]))$ are
$H_v$-submodules  of $M_2$. But
$U(f(\widetilde{M}_A);[t,s])=f(U(\widetilde{M}_A;[t,s]))$ and
$L(f(\widetilde{N}_A);[t,s]) =f(L(\widetilde{N}_A;[t,s]))$, so
$U(f(\widetilde{M}_A);[t,s])$ and $L(f(\widetilde{N}_A);[t,s])$
are $H_v$-submodules of $M_2$. Therefore $f[A]$ is an interval
valued intuitionistic $(S,T)$-fuzzy $H_v$-submodule of $M_2$.

\medskip
$(ii)$ For any $x,y\in H$ and $\alpha\in x+y$, we have
\[
\widetilde{M}_{f^{-1}(B)}(\alpha)=\widetilde{M}_B(f(\alpha))
 \ge T(\widetilde{M}_B(f(x)),\widetilde{M}_B(f(y)))
 =T(\widetilde{M}_{f^{-1}(B)}(x),\widetilde{M}_{f^{-1}(B)}(y)).
\]
Therefore
\[
\inf\limits_{\alpha\in x+y}\widetilde{M}_{f^{-1}(B)}(\alpha)\ge
T(\widetilde{M}_{f^{-1}(B)}(x),\widetilde{M}_{f^{-1}(B)}(y)).
 \]
For $x,a\in M_2$ there exists $y\in M_2$ such that $x\in a+y$.
Thus $f(x)\in f(a)+f(y)$ and
 \[
 T(\widetilde{M}_{f^{-1}(B)}(x),\widetilde{M}_{f^{-1}(B)}(a))
=T(\widetilde{M_B}(f(x)),\widetilde{M_B}(f(a)))\le
\widetilde{M_B}(f(y)) =\widetilde{M}_{f^{-1}(B)}(y).
\]
In the same manner, we can show that for $x,a\in M_2$ there exists
$z\in M_2$ such that $x\in z+a$ and
$T(\widetilde{M}_{f^{-1}(B)}(x),\widetilde{M}_{f^{-1}(B)}(a))\le
\widetilde{M}_{f^{-1}(B)}(z)$.

It is not difficult to see that, for all $x\in M_2$, $r\in R$ and
$\alpha\in r\cdot x$, we have
$\widetilde{M}_{f^{-1}(B)}(\alpha)=\widetilde{M} ( f^{-1}(
\alpha))\ge \widetilde{M} (f(x))=\widetilde{M}_{f^{-1}(B)}(x)$,
whence $\inf\limits_{\alpha\in r\cdot x}
\widetilde{M}_{f^{-1}(B)}(\alpha)\ge
\widetilde{M}_{f^{-1}(B)}(x)$. This completes the proof that
$\widetilde{M}_{f^{-1}(B)}$ is an interval valued $T$-submodule of
$M_1$.

Similarly, we can prove $\widetilde{N}_{f^{-1}(B)}$ is an interval
valued $S$-fuzzy $H_v$-submodule of $M_1.$  Therefore $f^{-1}(B)$
is an interval valued intuitionistic $(S,T)$-fuzzy $H_v$-submodule
of $M_1$. \hfill $\Box$

\bigskip

The mail tools in the theory of $H_v$-structures are the
fundamental relations. Consider an $H_v$-module $M$ over an
$H_v$-ring $R$. If the relation $\gamma^*$ is the smallest
equivalence relation on $R$ such that the quotient $R/\gamma^*$ is
a ring, we say that $\gamma^*$ is the {\it fundamental equivalence
relation} on $R$ and $R/\gamma^*$ is the {\it fundamental ring}.
The fundamental relation $\varepsilon^*$ on $M$ over $R$ is the
smallest equivalence relation on $M$ such that $ M/\varepsilon^*$
is a module over the ring $R/\gamma^*$ (see \cite{12, 7}).

Let $\mathcal{U}$ be the set of all expressions consisting of
finite hyperoperations of either on $R$ and $M$ or the external
hyperoperation applied on finite sets of $R$ and $M$. Then a
relation $\varepsilon$ can be defined on $M$ whose transitive
closure is the fundamental relation $\varepsilon^*$. The relation
$\varepsilon$ is as follows:
 \[
x\varepsilon y\Longleftrightarrow\{x,y\}\subseteq u \ \ \mbox{ for
some } \ \ u\in \mathcal{U}.
\]
Let us denote $\widetilde{\varepsilon}$ the transitive closure of
$\varepsilon$. Then we can rewrite the definition of
$\widetilde{\varepsilon}$ on $M$ as follows:
\[
 a\widetilde{\varepsilon} b\Longleftrightarrow \left\{\begin{array}{lll}
 \mbox{there exist } z_1,z_2,\ldots,z_{n+1}\in M \mbox{ and } u_1,u_2,\ldots,u_n\in \mathcal{U}\\[2pt]
 \mbox{such that } z_1=a,\; z_{n+1}=b,\; \mbox{ and }
 \{z_i,z_{i+1}\}\subseteq u_i \;\mbox{ for all } \ i=1,\ldots,n.
\end{array}\right.
 \]
The fundamental relation $\varepsilon^*$ is the transitive closure
of the relation $\varepsilon$ (see \cite{11}).

Suppose $\gamma^*(r)$ is the equivalence class containing $r\in R$
and $\varepsilon^*(x)$ is the equivalence class containing $x\in
M$. On $M/\varepsilon^*$, the sum $\oplus$ and the external
product $\odot$ using the $\gamma^*$ classes in $R$, are defined
as follows:
 \[
\varepsilon^*(x)\oplus\varepsilon^*(y)=\varepsilon^*(c) \ \ \mbox{
for all } c\in\varepsilon^*(x)+\varepsilon^*(y),
 \]
 \[
\gamma^*(r)\odot\varepsilon^*(x)=\varepsilon^*(d) \ \ \mbox{ for
all } d\in \gamma^*(r)\cdot \varepsilon^*(x).
\]
The kernel of the canonical map $\varphi: M\rightarrow
M/\varepsilon^*$ is called the {\it core} of $M$ and is denoted by
$\omega_M$. Here we also denote by $\omega_M$ the zero element of
the group $(M/\varepsilon^*,\oplus)$. Also, we have $\omega_M
=\varepsilon^*(0)$ and $\varepsilon^*(-x)=-\varepsilon^*(x)$ for
all $x\in M.$

\paragraph{Definition 3.7.}
Let $A=(\widetilde{M}_A,\widetilde{N}_A)$ be an interval valued
intuitionistic fuzzy set. The intuitionistic fuzzy set
$A/\varepsilon^*=(
\widetilde{M}_{\varepsilon*},\widetilde{N}_{\varepsilon*})$ is
defined as the pair of maps
 \[\left\{\begin{array}{c}
\widetilde{M}_{\varepsilon*}: M/\varepsilon^*\rightarrow
D[0,1], \\[3pt]
\widetilde{N}_{\varepsilon*}: N/\varepsilon^*\rightarrow D[0,1]
\end{array}\right.
\]
such that
\[
\widetilde{M}_{\varepsilon^*}(\varepsilon^*(x)) =
\left\{\begin{array}{clll}
\sup\limits_{a\in\varepsilon^*(x)}\widetilde{M}_A(a)
&& \mbox{ if \ $\varepsilon^*(x)\ne\omega_M$}\\[10pt]
[1,1] && \mbox{ otherwise}\end{array}\right.
\]
and
 \[
\widetilde{N}_{\varepsilon^*}(\varepsilon^*(x)) =
\left\{\begin{array}{cll}
\inf\limits_{a\in\varepsilon^*(x)}\widetilde{N}_A(a)
&& \mbox{ if \ $ \varepsilon^*(x)\ne\omega_M$}\\[10pt]
[0,0]&& \mbox{ otherwise}.\end{array}\right.
\]

\paragraph{Definition 3.8.}
Let $T$ (resp. $S$) be an idempotent interval $t$-norm (resp.
$s$-norm). An interval valued intuitionistic fuzzy set
$A=(\widetilde{M}_A,\widetilde{N}_A)$ on an ordinary module $M$
over a ring $R$ is called an {\it interval valued intuitionistic
$(S,T)$-fuzzy submodule} of $M,$ if
 \begin{enumerate}
\item[$(i)$] \ $\widetilde{M}_A(0)=[1,1]$ \ and \
$\widetilde{N}_A(0)=[0,0]$,
\item[$(ii)$] \
$T(\widetilde{M}_A(x),\widetilde{M}_A(y))\le\widetilde{M}_A(x-y)$
and
$S(\widetilde{N}_A(x),\widetilde{N}_A(y))\ge\widetilde{N}_A(x-y)$
for all $ x,y\in M,$
\item[$(iii)$] \ $\widetilde{M}_A(x)\le\widetilde{M}_A(r\cdot x)$
and $\widetilde{N}_A(x)\ge\widetilde{N}_A(r\cdot x)$ for all $x\in
M$ and $r\in R$.
\end{enumerate}

\paragraph{Theorem 3.9.} {\em
Let $M$ be an $H_v$-submodule of $M$ over an $H_v$-ring. If
$A=(\widetilde{M}_A,\widetilde{N}_A)$  is an interval valued
intuitionistic $(S,T)$-fuzzy $H_v$-submodule of $M$, then
$A/\varepsilon^* $  is an interval valued intuitionistic
$(S,T)$-fuzzy submodule of the fundamental module
 $M/\varepsilon^*$.}

\medskip\noindent
{\it Proof.} The first condition of the above definition is
trivially satisfied. To prove the second consider two arbitrary
elements $\varepsilon^*(x)$, $\varepsilon^*(y)$ of
$M/\varepsilon^*$.

If $\varepsilon^*(x)=\omega_M$, then
\[
T(\widetilde{M}_{\varepsilon^*}(\varepsilon^*(x)),
\widetilde{M}_{\varepsilon^*}(\varepsilon^*(y))\}=
T([1,1],\widetilde{M}_{\varepsilon^*}(\varepsilon^*(y)))
=\widetilde{M}_{\varepsilon^*}(\varepsilon^*(y))
=\widetilde{M}_{\varepsilon^*}(\varepsilon^*(x)\oplus\varepsilon^*(y)).
 \]

If $\varepsilon^*(x)\ne\omega_M$, then
\[\arraycolsep=.5mm\begin{array}{rl}
T(\widetilde{M}_{\varepsilon^*}(\varepsilon^*(x)),
\widetilde{M}_{\varepsilon^*}(\varepsilon^*(y)))&
=T(\sup\limits_{a\in\varepsilon^*(x)}\widetilde{M}_A(a) ,
\sup\limits_{b\in \varepsilon^*(y)}\widetilde{M}_A (b))
=\sup\limits_{\substack{a\in\varepsilon^*(x)\\b\in
\varepsilon^*(y)}}T(\widetilde{M}_A(a),\widetilde{M}_A(b))\\[14pt]
&\le\sup\limits_{\substack{a\in \varepsilon^*(x)\\b\in
\varepsilon^*(y)}}\big(\inf\limits_{\alpha\in
a+b}\widetilde{M}_A(\alpha)\big)\le\sup\limits_{\substack{a\in
\varepsilon^*(x)\\b\in\varepsilon^*(y)}}\big(\sup\limits_{\alpha\in
a+b}\widetilde{M}_A(\alpha)\big)\\[14pt]
&\le\sup\limits_{\substack{a\in
\varepsilon^*(x)\\b\in\varepsilon^*(y)}}\big(\sup\limits_{\alpha\in
\varepsilon^*(a+b)}\widetilde{M}_A(\alpha)\big)
=\sup\limits_{\substack{a\in\varepsilon^*(x)\\b\in
\varepsilon^*(y)}}\big(\widetilde{M}_{\varepsilon^*}(\varepsilon^*(a+b))\big)
\\[16pt]
&=\widetilde{M}_{\varepsilon^*}(\varepsilon^*(a+b))
\end{array}
\]
for all $a\in\varepsilon^*(x)$ and $b\in\varepsilon^*(y)$. Hence
\[
\widetilde{M}_{\varepsilon^*}(\varepsilon^*(a+b))=
\widetilde{M}_{\varepsilon^*}(\varepsilon^*(x+y))
=\widetilde{M}_{\varepsilon^*}(\varepsilon^*(x)\oplus\varepsilon^*(y)).
\]
So,
 \begin{equation}\label{r1}
T(\widetilde{M}_{\varepsilon^*}(\varepsilon^*(x)),\widetilde{M}_{\varepsilon^*}(\varepsilon^*(y)))
\le\widetilde{M}_{\varepsilon^*}(\varepsilon^*(x)\oplus\varepsilon^*(y)).
\end{equation}
The proof of the inequality
\begin{equation}\label{r2}
S(\widetilde{N}_{\varepsilon^*}(\varepsilon^*(x)),
\widetilde{N}_{\varepsilon^*}(\varepsilon^*(y)))
\ge\widetilde{N}_{\varepsilon^*}(\varepsilon^*(x)\oplus\varepsilon^*(y))
\end{equation}
is similar.

Let $\varepsilon^*(x)$ and $\varepsilon^*(a)$ be two arbitrary
elements of $M/\varepsilon^*$. Because
$A=(\widetilde{M}_A,\widetilde{N}_A)$ is an interval valued
intuitionistic $ (S,T)$-fuzzy $H_v$-submodule of $M$, then for
every $t\in\varepsilon^*(a)$ and $s\in \varepsilon^*(x)$, there
exists $y_{t,s}\in M$ such that $t\in s+y_{t,s}$ and
$T(\widetilde{M}_A(t),\widetilde{M}_A(s))\le\widetilde{M}_A(y_{t,s})$.
From $t\in s+y_{t,s}$, it follows that
$\varepsilon^*(s)\oplus\varepsilon^*(y_{t,s})=\varepsilon^*(t)$,
i.e.
$\varepsilon^*(x)\oplus\varepsilon^*(y_{t,s})=\varepsilon^*(a)$.

If $\varepsilon^*(a)\ne\omega_M$, $\varepsilon^*(x)\ne\omega_M$,
then
 \[
 \arraycolsep=.5mm\begin{array}{rl}
T(\widetilde{M}_{\varepsilon^*}(\varepsilon^*(a)),
\widetilde{M}_{\varepsilon^*}(\varepsilon^*(x)))&
=T(\sup\limits_{t\in\varepsilon^*(a)}\widetilde{M}_A(t),
\sup\limits_{s\in\varepsilon^*(a)}\widetilde{M}_A(s))
=\sup\limits_{\substack{t\in\varepsilon^*(a)\\s\in\varepsilon^*(x)}}
T( \widetilde{M}_A(t),\widetilde{M}_A(s))\\[22pt]
&\le\sup\limits_{\substack{t\in\varepsilon^*(a)\\s\in\varepsilon^*(x)}}
\widetilde{M}_A(y_{t,s})
\le\sup\limits_{y\in\varepsilon^*(y_{t,s})}\widetilde{M}_A(y)
=\widetilde{M}_{\varepsilon^*}(\varepsilon^*(y_{t,s})),
\end{array}
 \]
i.e. \[ T(\widetilde{M}_{\varepsilon^*}(\varepsilon^*(a)),
\widetilde{M}_{\varepsilon^*}(\varepsilon^*(x)))\le
 \widetilde{M}_{\varepsilon^*}(\varepsilon^*(y_{t,s})).
\]

Similarly
 \[
S(\widetilde{N}_{\varepsilon^*}(\varepsilon^*(a)),
\widetilde{N}_{\varepsilon^*}(\varepsilon^*(x)))\ge
\widetilde{N}_{\varepsilon^*}(\varepsilon^*(y_{t,s})) .
 \]

If $\varepsilon^*(x)=\omega_M(x)$, then
$\varepsilon^*(a)=\varepsilon^*(y_{t,s})$. So,
\[
T(\widetilde{M}_{\varepsilon^*}(\varepsilon^*(a)),
\widetilde{M}_{\varepsilon^*}(\varepsilon^*(x)))\le
\widetilde{M}_{\varepsilon^*}(\varepsilon^*(a))
=\widetilde{M}_{\varepsilon^*}(\varepsilon^*(y_{t,s}))
\]
and
\[
S(\widetilde{N}_{\varepsilon^*}(\varepsilon^*(a)),
\widetilde{N}_{\varepsilon^*}(\varepsilon^*(x)))\ge
\widetilde{N}_{\varepsilon^*}(\varepsilon^*(a))
=\widetilde{N}_{\varepsilon^*}(\varepsilon^*(y_{t,s})).
\]

Therefore for all $\varepsilon^*(x),\varepsilon^*(a)\in
M/\varepsilon^*$, there exists $\varepsilon^*(y)\in
M/\varepsilon^*$ such that
$\varepsilon^*(x)=\varepsilon^*(a)\oplus\varepsilon^*(y)$ for
which
 \[T(\widetilde{M}_{\varepsilon^*}(\varepsilon^*(x)),
\widetilde{M}_{\varepsilon^*}(\varepsilon^*(a)))\le
\widetilde{M}_{\varepsilon^*}(\varepsilon^*(y))
 \]
and
 \[
S(\widetilde{N}_{\varepsilon^*}(\varepsilon^*(x)),
\widetilde{N}_{\varepsilon^*}(\varepsilon^*(a)))\ge
\widetilde{N}_{\varepsilon^*}(\varepsilon^*(y)).
 \]

From the above it follows that for all $\varepsilon^*(x)\in
M/\varepsilon^*$ we have
$\widetilde{M}_{\varepsilon^*}(\varepsilon^*(x))\le
\widetilde{M}_{\varepsilon^*}(-\varepsilon^*(x))$ and
$\widetilde{N}_{\varepsilon^*}(\varepsilon^*(x))\ge
\widetilde{N}_{\varepsilon^*}(-\varepsilon^*(x))$. Indeed, for
$\omega_M\in M/\varepsilon^*$ there exists $\varepsilon^*(y_1)\in
M/\varepsilon^*$ such that
$\omega_M=\varepsilon^*(x)\oplus\varepsilon^*(y_1)$ and
\[
\widetilde{M}_{\varepsilon^*}(\varepsilon^*(x))=T([1,1],
\widetilde{M}_{\varepsilon^*}(\varepsilon^*(x))
=T(\widetilde{M}_{\varepsilon^*}(\omega_M),
\widetilde{M}_{\varepsilon^*}(\varepsilon^*(x)))
\le\widetilde{M}_{\varepsilon^*}(\varepsilon^*(y_1)),
\]
\[
\widetilde{N}_{\varepsilon^*}(\varepsilon^*(x))=
S([0,0],\widetilde{N}_{\varepsilon^*}(\varepsilon^*(x))=
S(\widetilde{N}_{\varepsilon^*}(\omega_M),
\widetilde{N}_{\varepsilon^*}(\varepsilon^*(x)))
\ge\widetilde{N}_{\varepsilon^*}(\varepsilon^*(y_1)), \] because
$\widetilde{M}_{\varepsilon^*}(\omega_M)=[1,1]$, \
$\widetilde{N}_{\varepsilon^*}(\omega_M)=[0,0]$. But
$\omega_M=\varepsilon^*(x)\oplus\varepsilon^*(y_1)$ implies
$\varepsilon^*(y_1)=-\varepsilon^*(x)$. Therefore
\begin{equation}\label{r3}
\widetilde{M}_{\varepsilon^*}(\varepsilon^*(x))\le
\widetilde{M}_{\varepsilon^*}(-\varepsilon^*(x)), \ \ \ \
\widetilde{N}_{\varepsilon^*}(\varepsilon^*(x))\ge
\widetilde{N}_{\varepsilon^*}(-\varepsilon^*(x)).
\end{equation}
So, by (\ref{r1}), (\ref{r2}) and (\ref{r3}), for all
$\varepsilon^*(x), \varepsilon^*(y)\in M/\varepsilon^*$, we have
 \[
 \arraycolsep=.5mm\begin{array}{rl}
\widetilde{M}_{\varepsilon^*}(\varepsilon^*(x)-\varepsilon^*(y))
&=\widetilde{M}_{\varepsilon^*}(\varepsilon^*(x)\oplus(-\varepsilon^*(y)))
\ge T(\widetilde{M}_{\varepsilon^*}(\varepsilon^*(x)),
\widetilde{M}_{\varepsilon^*}(-\varepsilon^*(y)))\\[4pt]
&\ge T(\widetilde{M}_{\varepsilon^*}(\varepsilon^*(x)),
\widetilde{M}_{\varepsilon^*}(\varepsilon^*(y)))
 \end{array}
 \]
and
 \[
 \arraycolsep=.5mm\begin{array}{rl}
\widetilde{N}_{\varepsilon^*}(\varepsilon^*(x)-\varepsilon^*(y))
&=\widetilde{N}_{\varepsilon^*}(\varepsilon^*(x)\oplus(-\varepsilon^*(y)))
\le S(\widetilde{N}_{\varepsilon^*}(\varepsilon^*(x)),
\widetilde{N}_{\varepsilon^*}(-\varepsilon^*(y)))\\[4pt]
&\le S(\widetilde{N}_{\varepsilon^*}(\varepsilon^*(x)),
\widetilde{N}_{\varepsilon^*}(\varepsilon^*(y))).
\end{array}
\]

This completes the proof of the second condition of Definition
3.8.

To prove the third condition observe that if $\varepsilon^*(x)\in
M/\varepsilon^*$ and $\gamma^*(r)\in R/\gamma^*$, then
$\widetilde{M}_{\varepsilon^*}(\gamma^*(r)\odot\varepsilon^*(x))=
\widetilde{M}_{\varepsilon^*}(\varepsilon^*(r\cdot b))$ for every
$b\in \varepsilon^*(x)$ and
\[
\widetilde{M}_{\varepsilon^*}(\varepsilon^*(r\cdot b))
=\sup\limits_{\alpha\in\varepsilon^*(r\cdot b)}
\widetilde{M}_A(\alpha)\ge\sup\limits_{\alpha\in r\cdot b}
\widetilde{M}_A(\alpha)\ge\widetilde{M}_A(b).
\]
Hence
$\widetilde{M}_{\varepsilon^*}(\gamma^*(r)\odot\varepsilon^*(x))
\ge\sup\limits_{b\in\varepsilon^*(x)}\widetilde{M}_A(b)$, which
implies
$\widetilde{M}_{\varepsilon^*}(\gamma^*(r)\odot\varepsilon^*(x))
\ge\widetilde{M}_{\varepsilon^*}(\varepsilon^*(x))$.

Similarly, we obtain
$\widetilde{N}_{\varepsilon^*}(\gamma^*(r)\odot\varepsilon^*(x))
\le\widetilde{N}_{\varepsilon^*}(\varepsilon^*(x))$. This
completes the proof. \hfill $\Box$

\end{document}